\mag=\magstep1 
\input epsf 
\documentclass[10pt,a4paper]{amsart} 
\usepackage{amssymb,latexsym} 

\def\Cal{\mathcal}

\def\<<{\langle } 
\def\>>{\rangle }

\numberwithin{equation}{section} 
\newtheorem{theorem}{Theorem}[section] 
\newtheorem{proposition}[theorem]{Proposition} 
\newtheorem{corollary}[theorem]{Corollary} 
\newtheorem{definition}[theorem]{Definition} 
 
\newtheorem{remark}[theorem]{Remark} 
\newtheorem{lemma}[theorem]{Lemma} 
\newtheorem{fact}[theorem]{Fact} 
\newtheorem{example}[theorem]{Example}

\newtheorem{notation}[theorem]{Notation}

\begin{document} 
 
\title{ 
Rational convex cones and cyclotomic multiple zeta values
} 
 
\author{Tomohide Terasoma} 

\date{\today} 
%
 
\maketitle 
 
\makeatletter 

\section{Introduction}

In this paper, we introduce zeta values of rational convex cones in
a finite dimensional vector space over $\bold Q$, which is a generalization
of cyclotomic multiple zeta values. Cyclotomic multiple zeta values appears in
the period integral for the fundamental group of $\bold C^{\times}-\mu_N$,
where $\mu_N$ is the subgroup of $N$-th root of unities in $\bold C^{\times}$.
The reader may find references for cyclotomic 
multiple zeta values in \cite{R}.

Let $C$ be a rational convex cone in $\bold R^m$ and $l_1, \dots, l_n$
be rational linear forms on $\bold Q^m$ which is positive
on the interior of $C$. Let $\chi$ be a finite order
character of $\bold Z^m$. The summation of
$$
\frac{\chi(x)}{l_1(x)\cdots l_n(x)}
$$
over the integral points $x$ in the interior of $C$ is denoted by
$\zeta_C(l_1, \dots, l_n,\chi)$ and called the zeta value of
the convex cone $C$. Our main theorem asserts that the zeta value
of a convex cone $C$ can be expressed as a $\bold Q^{ab}$-linear
combination of cyclotomic multiple zeta values,
where $\bold Q^{ab}$ denotes the extension of $\bold Q$
adjoining all roots of unities.
Here we give the outline of this paper. After defining zeta values
of convex cones, we show that the zeta values have integral expressions.
Motivated by these integral expressions, we introduce one variable functions
$I(y_n)$ and we compare them with cyclotomic multiple polylogarithms, whose
special values are cyclotomic multiple zeta values.
The zeta values of $C$ can be expressed as an integral of
$I(y_n)$'s. 
We study sequence of differential equations satisfied by
these functions $I(y_n)$'s and compare cyclotomic multiple polylogarithm
and $I(y_n)$'s.
Main theorem is a consequence of these relations
of functions.

By integral expressions of zeta values of cones, 
these zeta values can be regarded as period integrals of varieties
of pairs of some open parts of toric varieties and
sub-tori. Related varieties are defined as follows.
Let $(\bold C^{\times})^n$ be a torus and 
$$
\phi_i,\psi_j:(\bold C^{\times})^n \to \bold C^{\times}
$$
be surjective homomorphisms for $i=1, \dots, m, j=1,\dots, m'$.
We define divisors $D_i=\phi^{-1}(1),B_j=\psi_j^{-1}(1)$.
We set $U=(\bold C^{\times})^n-\cup_iD_i$ and 
$B=U\cap(\cup_jB_j)$.
Then the relative cohomology $H^i(U;B,\bold Q)$ is of mixed Tate type.
By our main theorem, the period integral of $H^i(U;B,\bold Q)$
can be expressed in terms of cyclotomic multiple zeta values.
In this paper, we produce concrete algorithm to get the expression 
of zeta values of cones by
cyclotomic multiple zeta values.

\section{Cyclotomic multiple zeta values and zeta values of convex cones} 

In this and next sections, we give the definition 
of zeta values of convex cones
and their integral expression.
\begin{definition}
\begin{enumerate}
\item
Let $N,k_1, \dots ,k_m \geq 1$ be
natural numbers and $\chi=(\chi_1, \dots ,\chi_m)$ be a homomorphism 
from $\bold Z^m$
to $\mu_N$. We assume that $k_m\neq 1 $ if $\chi (0, \dots,0, 1)$ is 
a trivial character. We define a cyclotomic multiple zeta value
$\zeta (k_1, \dots ,k_m, \chi)$ by
$$
\zeta (k_1,\dots ,k_m, \chi)
=
\sum_{a=(a_1, \dots, a_m) \in (\bold N^\times)^m}
\frac{\chi(a)}{a_1^{k_1} (a_1+a_2)^{k_2}
\cdots (a_1+\cdots +a_m)^{k_m}},
$$
where $\bold N^\times=\bold N-\{0\}$. The natural number $N$ is called
the modulus of the cyclotomic multiple zeta values.
\item
Let $\mu_N$ be the subgroup of $N$-th root of
unities in $\bold C^{\times}$ and $\mu_{\infty}=\cup_N\mu_N$.
We put $\bold Q^{ab} = \cup_{N \geq 1}\bold Q(\mu_N)$.
Let $\Cal Z_N$ be the $\bold Q^{ab}$-sub linear space 
of $\bold C$ generated by
multiple zeta values of modulus $N$, i.e.
\begin{align*}
\Cal Z_N=\langle
\zeta(k_1, \dots ,k_m, \chi) \mid &
\chi \text{ is a character whose values are contained in }\\
&
\mu_N, k_m \neq 1
\text{ if } \chi(0, \dots,0 ,1) =1
\rangle
\end{align*}
We define $\Cal Z_{\infty}$ by the linear hull of $\Cal Z_N$ for all
$N \geq 1$.
\end{enumerate}
\end{definition}
\begin{definition}[Rational convex cone]
\begin{enumerate}
\item
A closed convex set $C$ in $\bold R^m$ is called a rational convex cone
if there exist finite number of rational vectors $v_1, \dots ,v_b$
such that 
$$
C=\bold R_{+}v_1 + \cdots + \bold R_{+}v_b.
$$
Moreover, if we can choose linearly independent vectors 
$v_1, \dots, v_b$ over $\bold R$,
$C$ is called a simplicial cone.
The interior of the cone $C$ in the $\bold R$-linear hull
of $C$ is denoted by $C^0$. 
\item
A subset $\sigma$ of a convex cone $C$ is called a 
face of $C$ if there exists a linear form $l$ of $\bold Q^m$
such that (1) $l (\sigma ) \geq 0$, and (2) 
$\sigma = \{x \mid l(x)=0\}\cap C$.
\item
Let $\chi$ be a
character of $\bold Z^m$ and $l_1, \dots ,l_n$ be rational 
$\bold Q$-linear forms on $\bold Q^m$ such that $l_i(C^0) >0$
($i=1, \dots ,n$). We define $\zeta_C(l_1, \dots ,l_n, \chi)$
by
$$
\zeta_C(l_1, \dots , l_n, \chi) 
=
\sum_{x 
\in C^0\cap \bold Z^m}\frac{\chi(x)}
{l_1(x)\cdots l_n(x)},
$$
if the sum is absolutely convergent.
$\zeta_C(l_1, \dots ,l_m, \chi)$ is called a zeta value of
a convex cone $C$.
\end{enumerate}
\end{definition}
The main theorem of this paper is the following.
\begin{theorem}
The value $\zeta_C(l_1, \dots ,l_n, \chi)$
is an element of $\Cal Z_{\infty}$.
\end{theorem}
\begin{remark}
\begin{enumerate}
\item
Cyclotomic multiple zeta values that converge absolutely
are special case of multiple zeta values
for convex cones and linear forms.
\item
A proto-type of this summation can be found in Zagier's paper \cite{Z}.
\item
In this paper, we give an explicit algorithm to compute
the value of $\zeta_C(l_1, \dots, l_n,\chi)$.
\end{enumerate}
\end{remark}

\section{First reduction and integral expression}

In this section, we reduce the main theorem to the case where $C$
is a simplicial cone and the semi group $C\cap \bold Z^m$ is freely
generated by integral $m$ elements. Moreover we show that
the zeta value of a convex cone $C$ admits an integral representation.

\begin{lemma}
\label{induced rep}
Let $L \supset \bold Z^m$ be a lattice in $\bold Q^m$ such that 
the index $\kappa=[L:\bold Z^m]$ 
is finite. Let $C$ be a rational convex cone in
$\bold R^m$. Let $\chi :\bold Z^m \to \mu_N$ be a character and
$$
Ind_{\bold Z^m}^L(\chi) =\oplus _{i=1}^{\kappa}\chi_i
$$
be the absolutely irreducible decomposition of $Ind_{\bold Z^m}^L(\chi)$.
Then we have
$$
\sum_{x\in C^{0} \cap \bold Z^m}\frac{\chi (x)}{l_1(x)\cdots l_n(x)}
=
\frac{1}{\kappa}\sum_{i=1}^{\kappa}
\sum_{x\in C^{0} \cap L}\frac{\chi_i(x)}{l_1(x)\cdots l_n(x)}.
$$
Moreover if the left hand side converges absolutely, then
each sums of the right hand side converges absolutely.
\end{lemma}
\begin{proof}
For $x\in L$, we have an equality
$$
Ind_{\bold Z^m}^{L}(\chi)(x)=\begin{cases}
[L:\bold Z^m]\chi(x) &\text{ if } x\in \bold Z^m \\
0 & \text{ if } x \notin \bold Z^m.
\end{cases}
$$
Let $x=(x_1, \dots, x_m)$ be a coordinate of $\bold Z^m$.
The condition for absolutely convergence of the summation 
for the left hand side (resp. right hand side) is equivalent to
the finiteness of the integral
$$
\int_D
\frac{1}{l_1(x)\cdots l_n(x)}dx_1\cdots dx_m
$$ 
where $D$ is the convex hull of 
$C^0\cap \bold Z^m$ (resp. $C^0 \cap L$).
Therefore two conditions are equivalent.
\end{proof}
\begin{proposition}
\label{integral express}
Let $C=\bold R_{+}^m$. Let $l_1, \dots, l_n$ be non-zero linear forms on 
$\bold Q^m$ defined by
\begin{align*}
l_1 &=l_{11}x_1 + \dots + l_{1m}x_m \\
 & \dots \\
l_n &=l_{n1}x_1 + \dots + l_{nm}x_m.
\end{align*}
with $l_{ij} \in \bold N$.
Let $\chi$ be a character of
$\bold Z^n$ of finite order. We put 
$$
M_1=y_1^{l_{11}}\cdots y_n^{l_{n1}},\dots ,
M_m=y_1^{l_{1m}}\cdots y_n^{l_{nm}}.
$$
Then we have
\begin{equation}
\label{integral expression 1}
\zeta_C(l_1, \dots ,l_n, \chi)
=
\int_{(0,1)^n}
 \frac{\chi (1,\dots ,1)M_1\cdots M_m}
{(1-\chi(\bold e_1)M_1)\cdots (1-\chi (\bold e_m)M_m)}
\frac{dy_1}{y_1}\cdots \frac{dy_n}{y_n}
,
\end{equation}
where $\bold e_i=(0,\dots, 0,\overset{\overset{i}\smile}{1},0,\dots ,0)$,
if the left hand side exists.
\end{proposition}
\begin{proof}
For $x \in \bold Z^m$, we have
\begin{align*}
\frac{\chi (x)}{l_1(x)\cdots l_n(x)} 
& =\int_{(0,1)^n}
\chi(x)
y_1^{l_1(x)} \cdots
y_n^{l_n(x)}
\frac{dy_1\cdots dy_n}{y_1\cdots y_n} \\
&=
\int_{(0,1)^n}
\chi (x)
M_1^{x_1} \cdots
M_m^{x_m}
\frac{dy_1\cdots dy_n}{y_1\cdots y_n}.
\end{align*}
By taking the summation for $x\in (\bold N^{\times})^m$
and using the absolutely
convergent condition for summation, we change the integration and
summation. Thus we get the theorem.
\end{proof}

\begin{corollary}
For any rational cone $C$, character $\chi$ of $\bold Z^m$ and
linear forms $l_1, \dots ,l_n$ with $l_i(C^0) >0$, the zeta value
of $\zeta_C(l_1, \dots, l_n, \chi)$
can be expressed as a $\bold Q^{ab}$-linear combination of integrals of
the form (\ref{integral expression 1}).
\end{corollary}
\begin{proof}
We choose a decomposition $C=\cup_{i=1}^k C_i$ of rational convex cone $C$, 
by rational $m$-dimensional simplicial cones $C_i$. Then $C^0$ is a disjoint
union of open part of simplicial cones, which is denoted as
$$
C^0 =\coprod_{i=1}^{\kappa}D_i
\quad (D_i\simeq {\bold R^{\times}_{+}}^{d_i}).
$$
Let $V_i$ be a linear hull of $D_i$ and $V_i\cap \bold Z^m$ is denoted
as $L_i$. The zeta values of $D_i$ with respect to the lattice 
$L_i$ is denoted as
$\zeta_{D_i}(l_1\mid_{D_i}, \dots, l_n\mid_{D_i},\chi\mid_{D_i})$.
Then we have
\begin{equation}
\label{simplicial decomp}
\zeta_C(l_1, \dots, l_n,\chi) =
\sum_{i=1}^{\kappa}
\zeta_{D_i}(l_1\mid_{D_i}, \dots, l_n\mid_{D_i},\chi\mid_{D_i}).
\end{equation}
For each $i$, we can choose a super-lattice $\tilde L_i \supset L_i$
with $[\tilde L_i:L_i] < \infty$ and $\tilde L_i\cap \overline{D_i}$
is generated by $d_i$ elements. By using  Lemma \ref{induced rep}
and Proposition \ref{integral express},
each term of the right hand side of
(\ref{simplicial decomp}) can be express as a $\bold Q^{ab}$-linear
combination of the integrals of the form (\ref{integral expression 1}).
\end{proof}
\begin{remark}
Harmonic shuffle relations for multiple zeta values come from
these decompositions.
\end{remark}
\section{Preparation for a decomposition by simplicial cones}

In this section, we give some construction of
decomposition of a convex cone by rational simplicial cones.
Let $\alpha_1, \dots , \alpha_l$  be non-zero linear forms on $\bold Q^n$
and $C$ an $n$-dimensional rational convex cone.
\begin{notation}
Let $\alpha$ be a linear form on $\bold Q^n$. The set 
$\alpha\cdot \bold Q^{\times}$ is called a linear form up to constant
multiple and denoted as $[\alpha]$.
\end{notation}
\begin{definition}
\begin{enumerate}
\item
Let $\Delta$ be a simplicial cone and $\alpha$ be a linear form on $\Delta$.
The linear form $\alpha$ is said to be definite on $\Delta$ if
$\alpha (\Delta) \geq 0$ or $\alpha (\Delta) \leq 0$.
\item
Let $S=\{\bar\alpha_1, \dots ,\bar\alpha_l\}$ be a set of linear forms
up to constant multiple.
A decomposition 
\begin{equation}
\label{decomposition by simplicials}
C= \cup_{i=1}^k\Delta_i
\end{equation}
of $C$ by simplicial cones is said to be compatible
with respect to the set $S$
if $\bar\alpha_j$ is definite on $\Delta_i$ for $j=1, \dots, l$.
\item
Let $\alpha_1, \dots ,\alpha_l$ be linear forms on $\Delta$
and $\bar\alpha_i=[\alpha_i]$.
Let $F$ be a codimension one face of
$\Delta$. The set $S=\{\bar\alpha_1, \dots ,\bar\alpha_l\}$ 
of linear forms up to constant is 
said to be non-degenerate with respect
to $F$ if the restriction of $\alpha_i$ to $F$ is non-zero for $i=1, \dots ,l$.
\end{enumerate}
\end{definition}
We use the following fact
\begin{fact}
Let $C$ and $S=\{\bar\alpha_1, \dots ,\bar\alpha_l\}$ 
be as above. Then there exists
a decomposition (\ref{decomposition by simplicials})
with the following properties.
\begin{enumerate}
\item
The decomposition is compatible with respect to $S$.
\item
There exists a codimension one face $F_i$ of $\Delta_i$
such that the set $S=\{\bar\alpha_1, \dots ,\bar\alpha_l\}$ is
non-degenerate with respect to $F_i$.
\end{enumerate}
\end{fact}
\begin{definition}
\begin{enumerate}
\item
A sequence of simplicial cones 
$$
F:\Delta=\Delta^{(0)}\supset
\Delta^{(1)} \supset \cdots \supset \Delta^{(n-1)}
$$
is called a flag of $\Delta$ if $\Delta^{(i)}$ is a
codimension $i$ face of $\Delta^{(i)}$.
For a flag $F$, there exists a coordinate $(\eta_1, \dots, \eta_n)$
of $\bold Q^n$ such that
\begin{enumerate}
\item
$
\Delta =\{\eta_i \geq 0 \text{ for }i=1, \dots, n\},
$
\item
$
\Delta^{(i)}=\{(\eta_1, \dots ,\eta_n)\in \Delta \mid 
\eta_1= \cdots =\eta_i=0 \}.
$
\end{enumerate}
This coordinate is called a standard coordinate for the flag $F$.
\item
Let $\Delta$ be an $n$-dimensional simplicial cone, $\sigma$
be an $(n-i)$-dimensional face of $\Delta$. Then there exist a unique
$i$ dimensional face $\sigma^*$ of $\Delta$ 
such that $\sigma \cap \sigma^*=\{0\}$.
The face $\sigma^*$ is called the dual face of $\sigma$.
\item For two convex cones $\sigma_1$ and $\sigma_2$,
the linear join $\sigma_1 * \sigma_2$ 
of $\sigma_1$ and $\sigma_2$ is defined by
$$
\sigma_1 * \sigma_2=
\{ ax+by\mid a, b\in \bold R, a, b  \geq 0, a+b=1, x \in \sigma_1,
y \in \sigma_2\}.
$$
\end{enumerate}
\end{definition}
\begin{remark}
\begin{enumerate}
\label{rem for face and index}
\item
If two linear forms $\alpha_1$ and $\alpha_2$ 
are distinct up to constant multiple
and $\alpha_1(v) \neq 0, \alpha_2(v) \neq 0$,
Then the restriction 
$(\alpha_1(v)\alpha_2-\alpha_2(v)\alpha_1)\mid_F$ is non-zero.
\item
Let $F$ be a flag of $n$-dimensional
simplicial cone $\Delta$ and $\Delta^{(n-1)}$
the one dimensional component of the
flag $F$. The dual of $\Delta^{(n-1)}$
in $\Delta$ is a codimension one face and denoted as $ir(F)$.
The face $ir(F)$ is called the irregular face of the flag $F$.
A face $\sigma$ of $\Delta$ is called a regular face
if it is not contained in the irregular face $ir(F)$.
Let $(\eta_1, \dots,\eta_n)$ be a standard coordinate for the flag
$F$. Then an $n_{\sigma}$-dimensional regular face $\sigma$
is defined by 
\begin{equation}
\label{equation for sigma}
\eta_j =0
\text{ for }
j \notin\{i_1, \dots ,i_{n_{\sigma}}\},
\end{equation}
where $1\leq i_1< \dots < i_{n_{\sigma}}=n$.
The coordinate $(\eta_{i_1}, \dots, \eta_{i_{n_\sigma}})$ of $\sigma$
is called the standard coordinate.
We put $\eta_{i_1}=\eta_1^{\sigma},\cdots,\eta_{i_{n_{\sigma}}}
=\eta_{n_{\sigma}}^{\sigma}$.
\end{enumerate}
\end{remark}

We define the notion of derived sequence of linear forms.
\begin{definition}[Derived sequence of linear forms]
\label{def derived seq}
\begin{enumerate}
\item
\label{derived set}
Let $\Delta$ be a simplicial cone, $F$ be a codimension one face of $\Delta$
and $S=\{\bar\alpha_1, \dots ,\bar\alpha_l\}$ be a set of 
linear forms up to constant multiple
on $\Delta$ which is non-degenerate with respect to $F$.
Let $v$ be a non-zero vector in the dual $F^*$ of $F$.
We define a set $D_F(S)$ of linear forms on $F$ up to constant multiple
as the union of 
\begin{enumerate}
\item
the set $\{[\alpha\mid_F] \mid [\alpha] \in S\}$, and
\item
the set $\{[(\alpha_1(v)\alpha_2 - \alpha_2(v)\alpha_1)
\mid_F] \mid [\alpha_1], [\alpha_2] \in S,
\alpha_1(v) \neq 0, \alpha_2 \neq 0
\}$.
\end{enumerate}

\item
Let $F=\{\Delta^{(0)}\supset
\Delta^{(1)} \supset \cdots \supset \Delta^{(n-1)}\}$
be a flag of $\Delta$ and $S^{(i)}$ be a set of linear forms
up to constant multiple on $\Delta^{(i)}$. 
The pair $\{F, \{S^{(i)}\}_i\}$ is called
a derived sequence of linear forms if the following
conditions are satisfied.
\begin{enumerate}
\item
The set $S^{(i)}$ is non-degenerate with respect to the face 
$\Delta^{(i+1)}$.
\item
The set $S^{(i)}$ is definite on $\Delta^{(i)}$.
\item
$D_{\Delta^{(i+1)}}(S^{(i)}) \subset S^{(i+1)}$.
\end{enumerate}
\end{enumerate}
\end{definition}
\begin{proposition}
Let $F$ be a flag of $n$-dimensional simplicial cone $\Delta$
and $(F, \{S^{(i)}\}_i)$ a derived sequence of linear forms.
\begin{enumerate}
\item
The restriction of an element $\alpha \in S^{(i)}$ 
($i=1, \dots , n-1$) to
$\Delta^{(n-1)}$ is non-zero.
In other words, 
using a standard coordinate $(\eta_1, \dots, \eta_n)$ of
$\Delta$,
the linear form $\alpha$ can be written as
$$
\alpha=\alpha^{(i+1)}\eta_{i+1}+\dots +\alpha^{(n)}\eta_n
$$
with $\alpha^{(i+1)}\geq 0, \dots ,\alpha^{(n-1)} \geq 0, \alpha^{(n)} >0$.
\item
Let $\sigma$ be a regular face of the flag $F$.
The restriction of an element $\alpha \in S^{(i)}$ to 
$\sigma\cap \Delta^{(i)}$ is non-zero.
\item
Let $\sigma$ be a codimension one
regular face of the flag $F$ defined by $\eta_p=0$.
We define a flag of $F_{\sigma}=\sigma\cap F$ by
$$
\begin{matrix}
\sigma =\Delta^{(0)}\cap \sigma & \supset &
\Delta^{(1)}\cap \sigma  & \supset 
\cdots \supset &  \Delta^{(p-1)}\cap \sigma
 =  \Delta^{(p)} & \supset \cdots & \Delta^{(n-1)} \\
\parallel & &
\parallel  &  
 &  \parallel
 &  & \parallel \\
\Delta^{(0)}_{\sigma} & \supset &
\Delta^{(1)}_{\sigma}  & \supset 
\cdots \supset &  \Delta^{(p-1)}_{\sigma}
 & \supset \cdots & \Delta^{(n-2)}_{\sigma} \\
\end{matrix}
$$
For $i<p-1$, we define a set of linear forms $S^{(i)}_{\sigma}$ 
on $\Delta^{(i)}\cap\sigma$ by 
$$
S^{(i)}_{\sigma}=\{[\alpha\mid_{\Delta^{(i)}\cap \sigma}] 
\mid [\alpha] \in S^{(i)}\}.
$$
and for $i\geq p-1$, we define $S^{(i)}_{\sigma}=S^{(i+1)}$.
Then $(F_{\sigma},\{S^{(i)}_{\sigma}\}_{i=1, \dots ,n-2})$ is 
a derived sequence of linear forms.
\end{enumerate}
\end{proposition}
\begin{proof}
The statement (1) and (2) is obvious from the definition of derived
sequence of linear forms.

(3)
Let $\sigma^*$ be the dual face of $\sigma$. Then
$$
\Delta^{(i-1)}=\sigma^* *\Delta^{(i-1)}_\sigma \supset
\Delta^{(i)}=\sigma^* *\Delta^{(i)}_\sigma
$$
for $i\leq p-1$.
Therefore we can use common $v$ 
in 
Definition \ref{def derived seq} (\ref{derived set})
to define 
$D_{\Delta^{(i)}_{\sigma}}(S^{(i-1)}_{\sigma})$ and
$D_{\Delta^{(i)}}(S^{(i-1)})$.
As a consequence, for $i \leq p-1$, we have
\begin{align*}
D_{\Delta^{(i)}_{\sigma}}(S^{(i-1)}_{\sigma})
& =
\{ [\alpha \mid_{\Delta^{(i)}_{\sigma}}] \mid  [\alpha] \in
D_{\Delta^{(i)}}(S^{(i-1)})\} \\
& \subset
\{ [\alpha \mid_{\Delta^{(i)}_{\sigma}}] \mid  [\alpha] \in
S^{(i)}\}.
\end{align*}
Since
$\{ [\alpha \mid_{\Delta^{(i)}_{\sigma}}] \mid  [\alpha] \in
S^{(i)}\}=S^{(i)}_\sigma$ if $i<p-1$, and
$$
\{ [\alpha \mid_{\Delta^{(p-1)}_{\sigma}}] \mid  [\alpha] \in
S^{(p-1)}\}
\subset S^{(p)}=
S^{(p-1)}_{\sigma},
$$
the third condition for a derived sequence for linear forms are
satisfied.
The first condition is a consequence of the second statement of the
proposition, and
the second condition is obvious.
\end{proof}
\begin{corollary}
Let $\Delta$ be an $n$-dimensional simplicial cone, $F$ be
a flag of $\Delta$ and $\Cal D=(F, \{S^{(i)}\})$ be a derived sequence
for $F$.
Let $\sigma$ be an $m$-dimensional regular face for $F$ and 
$1 \leq i_1 < i_2 <\cdots < \eta_{i_{n_\sigma}}=n$ be the 
index set given in (\ref{equation for sigma}).
We set 
\begin{enumerate}
\item
$\Delta^{(0)}_{\sigma}=\sigma, 
\Delta^{(1)}_{\sigma}=\sigma\cap\Delta^{(i_2-1)},
\dots, \Delta^{({n_\sigma}-1)}_{\sigma}
=\sigma\cap\Delta^{(i_{n_\sigma}-1)}$ and
\item
$$
S^{(0)}_{\sigma}=S^{(i_1-1)}\mid_{\Delta^{(0)}_{\sigma}},
S^{(1)}_{\sigma}=S^{(i_2-1)}\mid_{\Delta^{(1)}_{\sigma}}, 
\dots,
S^{({n_\sigma}-1)}_{\sigma}
=S^{(i_{n_\sigma}-1)}\mid_{\Delta^{({n_\sigma}-1)}_{\sigma}}.
$$
\end{enumerate}
Then $F_{\sigma}=\{\Delta_{\sigma}^{(i)}\}_{i=0, \dots, n_\sigma-1}$
is a flag of $\sigma$ and 
$\Cal D_{\sigma}=(F_{\sigma},\{S_{\sigma}^{(i)}\}_{i=0, \dots, n_\sigma-1})$
is a derived sequence of $\sigma$.
\end{corollary}
\begin{definition}
The above derived sequence $\Cal D_{\sigma}$ is called the
restriction of $\Cal D$ to $\sigma$.
\end{definition}
\begin{proposition}
\label{preparation proposition}
Let $C$ be an $n$-dimensional rational convex cone in $\bold R^n$ and
$S=\{\bar\alpha_1, \dots ,\bar\alpha_l\}$ be a set of linear 
forms on $\bold Q^n$ distinct to each other
up to constant multiple.
Then there exists 
\begin{enumerate}
\item
a decomposition
$C=\cup_i^{k}\tilde\Delta_i$ by rational simplicial cones and
\item
a flag
$$
F_i:\tilde\Delta_i=\tilde\Delta^{(0)}_i\supset 
\tilde\Delta^{(1)}_i \supset \cdots \supset
\tilde\Delta^{(n-1)}_i
$$
for each $i=1, \dots ,k$, and
\item
a sequence of sets of linear forms
$S_i^{(0)},S_i^{(1)},\dots ,S_i^{(n-1)}$ up to constant multiple on
$\tilde\Delta^{(0)}_i,\tilde\Delta^{(1)}_i,\dots ,\tilde\Delta^{(n-1)}_i$
\end{enumerate}
such that pairs $\{F_i, \{S_i^{(j)}\}_j\}$ is a derived sequence
for $i=1, \dots ,k$ and $S_i^{(0)}=\{\bar\alpha_1, \dots ,\bar\alpha_l\}$.
\end{proposition}
\begin{proof}
Step 1.
We inductively construct the following sets and maps. 
\begin{enumerate}
\item
Index sets $I^{(0)}, I^{(1)}, \dots ,I^{(n-1)}$
and a surjective map $\rho:I^{(i)}\to I^{(i-1)}$ for $i=1,\dots ,n-1$. 
\item
$(n-i)$-dimensional simplicial cones $\Delta_j^{(i)}$ 
indexed by $j \in I^{(i)}$.
\item
A set $S_j^{(i)}$ of linear forms up to constant multiple
on $\Delta_j^{(i)}$.
\end{enumerate}
with the following properties:
\begin{enumerate}
\item
For $k\in I^{(i-1)}$,
$$
E_k^{(i-1)}=\cup_{\rho(j)=k}\Delta_j^{(i)}
$$ 
is a decomposition by simplicial cones
of a codimension one face of $\Delta_{k}^{(i-1)}$.
\item
The set $S_{k}^{(i-1)}$ is non-degenerate with respect
to the face $E_k^{(i-1)}$ of $\Delta_k^{(i-1)}$.
\item
Each element in $S_k^{(i-1)}$ is definite on $\Delta_k^{(i-1)}$.
\item
$S_j^{(i)}$ is equal to $D_{E_{\rho(j)}^{(i-1)}}(S_{\rho (j)}^{(i-1)})$.
\end{enumerate}

\noindent
Step 2.
Let $\rho_i: I^{(n-1)} \to I^{(i)}$ be the successive composite of
$\rho$.
We inductively construct the following decomposition of 
$\Delta^{(i)}_j$ for $j\in I^{(i)}$:
\begin{equation}
\label{inductive decomposition}
\Delta^{(i)}_j=\cup_{\rho_i (p)=j}\tilde\Delta^{(i)}_p.
\end{equation}
For $i=n-1$, we put $\tilde\Delta^{(n-1)}_p=\Delta^{(n-1)}_p$ 
for $p\in I^{(n-1)}$.
We construct a decomposition of $\Delta^{(i-1)}_k$ for $k\in I^{(i-1)}$
using decompositions (\ref{inductive decomposition})
of $\Delta^{(i)}_j$ for $j\in I^{(i)}$. The dual simplex of
$E^{(i-1)}_k$ in $\Delta^{(i-1)}_k$ is denoted by $G^{(i-1)}_k$.
For $p\in I^{(n-1)}$, we put 
$\tilde \Delta^{(i-1)}_p=\tilde \Delta^{(i)}_p*G^{(i-1)}_{\rho_{i-1}(p)}$.
Then we have
\begin{align*}
\Delta^{(i-1)}_k
& =E^{(i-1)}_k*G^{(i-1)}_k \\
& =(\cup_{\rho (j)=k}\Delta^{(i)}_j)*G^{(i-1)}_k \\
& =(\cup_{\rho_{i-1}(p)=k}\tilde\Delta^{(i)}_p)*G^{(i-1)}_k \\
& =\cup_{\rho_{i-1}(p)=k}\tilde\Delta^{(i-1)}_p.
\end{align*}
As a consequence, we have a decomposition
$$
\Delta^{(0)}=\cup_{p\in I^{(n-1)}}\tilde\Delta_p^{(0)}.
$$
of $\Delta^{(0)}$. On the simplicial cone $\tilde\Delta_p^{(i-1)}$,
we have a flag
$$
F_p:\tilde\Delta_p^{(0)}\supset \tilde\Delta_p^{(1)}\supset
\cdots \supset
\tilde\Delta_p^{(n-1)}.
$$
For $p \in I^{(n-1)}$,
we have a inclusion 
$\tilde \Delta_p^{(i)}\subset \Delta_{\rho_i(p)}^{(i)}$.
The restriction of $S^{(i)}_{\rho_i(p)}$ to $\tilde \Delta_p^{(i)}$
is denoted as $S^{(i)}_p$. Then we have
$$
S^{(i)}_p=D_{\tilde\Delta^{(i-1)}_p}(S^{(i-1)}_p)
$$
and the pair $(F_j,\{S^{(i)}_j\}_j)$ is a derived sequence of linear forms.
\end{proof}

Let $\Delta$ be an $n$-dimensional simplicial cone, $F$ be
a flag on $\Delta$ and $\Cal D=(F, \{S^{(i)}\}_i)$ be a derived sequence
of linear forms. For a regular face $\sigma$ of $\Delta$, the
restriction of $\Cal D$ to $\sigma$ is denoted as $\Cal D_{\sigma}$.
We fix a standard coordinate $(\eta_1, \dots ,\eta_n)$ of $\Delta$
with respect to the flag $F$. 
Let $e_1, \dots, e_n \geq 1$ be natural numbers.
Then the coordinate $(\tilde\eta_1,\dots, \tilde\eta_n)$
defined by $\eta_i=e_i\tilde\eta_i$ is also a standard coordinate for $F$.
Then the linear form 
$\alpha= \alpha^{(1)}\eta_1 + \cdots +\alpha^{(n)}\eta_n$
is transformed into
$$
\alpha= (e_1\alpha^{(1)})\tilde\eta_1 + \cdots +
(e_n\alpha^{(n)})\tilde\eta_n
$$
with respect to the coordinate $(\tilde\eta_1, \dots ,\tilde\eta_n)$.
\begin{lemma}
\label{changing variable}
Let $\Cal D=(F,\{S^{(i)}\})$ be a derived sequence of linear forms on 
$F$, $(\eta_1, \dots ,\eta_n)$ be a standard coordinate of $\Delta$
for $F$.
Then there exist natural numbers 
$e_1 , \dots ,e_n\geq 1$ with the following property:
By changing coordinate $(\tilde\eta_1, \dots, \tilde\eta_n)$
with a relation $\eta_i= e_i\tilde\eta_i$, any
element $\bar\alpha \in S^{(i)}_{\sigma}$
such that
$\bar\alpha(v) \neq 0$ for a non-zero vector 
$v$ in $(\Delta^{(i+1)}_{\sigma})^*$
has a representative 
$$
\alpha=\alpha_{i+1}\tilde\eta_{i+1}^{\sigma}+\cdots +
\alpha_{n_{\sigma}}\tilde\eta_{n_{\sigma}}^{\sigma}
$$
such that
$\alpha_{i+1}=1$ and $\alpha_j \in \bold N$ for $j \geq i+1$.
\end{lemma}
\begin{definition}
\begin{enumerate}
\item
The representative $\alpha=(\alpha_{i+1}, \dots, \alpha_n)$
with the above property is called the primitive representative.
The standard coordinate $(\eta_1, \dots, \eta_n)$ is called primitive
with respect to the derived sequence $\Cal D$ if it has the property of the
above lemma.
\item
Let $\Cal D=(F, \{S^{(i)}\}_i)$ be a derived sequence of linear forms
and $v$ a non-zero element of $(\Delta^{(i+1)})^*$. The set of primitive
representative of the subset
$$
\{\bar\alpha \in S^{(i)} \mid \bar\alpha (v)\neq 0\}
$$
of $S^{(i)}$
is called the variable part of $S^{(i)}$ and denoted as $S_{var}^{(i)}$.
\end{enumerate}
\end{definition}

\section{The second reduction, Definite integral of type $S$ 
and type $\Cal D$}

Let $\Delta$ be an $n$-dimensional simplicial cone and $F$
a flag of $\Delta$,
We choose a primitive standard coordinate $(\eta_1, \dots , \eta_n)$
of $\Delta$ with respect to the flag $F$.
Let $S=\{\bar\alpha_1, \dots, \bar\alpha_l\}$ be a set 
of distinct definite non-zero rational linear forms
on $\Delta$ up to constant.
For a linear form 
\begin{equation}
\label{integral form}
\alpha=\alpha^{(1)}\eta_1 + \dots + \alpha^{(n)}\eta_n
\text{ on $\Delta $ such that $\alpha^{(i)}\in \bold N$, }
\end{equation}
$y^{\alpha}$ denotes a monomial in $\bold Q[y_1, \dots,y_n]$
defined by
$y^{\alpha}=y_1^{\alpha^{(1)}}\cdots y_n^{\alpha^{(n)}}$.
\begin{definition}
\begin{enumerate}
\item
We consider a rational function of the form
\begin{equation}
\label{type of denom}
L=\prod_{i=1}^p
\Big(
\frac{e_iy^{\alpha_i}}{1-e_iy^{\alpha_i}}
\Big)
^{\mu_i},
\end{equation}
where 
$\alpha_i$ is a linear form on $\Delta$
with the condition (\ref{integral form})
and
$e_i$ is an element of $\mu_{N}$ for some $N\in \bold N$.
The integral
\begin{equation}
\label{definite integral}
I=\int_{(0,1)^n}L(y_1, \dots, y_n)\frac{dy_1}{y_1}\cdots \frac{dy_n}{y_n}
\end{equation}
is called a definite integral of type $S$ if 
$[\alpha_i] \in S$ ($i=1, \dots, p$)
and the integral
(\ref{definite integral}) converges absolutely.
\item
Let $F$ be a flag of $\Delta$.
For a derived sequence
$\Cal D=(F, \{S^{(i)}\}_i)$,
the integral (\ref{definite integral})
is called a definite integral of type $\Cal D$
if $\alpha_i \in\coprod_{i=1}^{n-1}S^{(i)}_{var}$
and the integral converges absolutely.
\end{enumerate}
\end{definition}
\begin{example}
\end{example}
The integral of the right hand side of (\ref{integral expression 1})
is a definite integral of type $S$, where
\begin{align*}
S &=\{[\alpha_1],\dots, [\alpha_m]\}, \\
\alpha_i &=l_{1i}\eta_1+\cdots + l_{ni}\eta_n\qquad (i=1, \dots, m).
\end{align*}

\begin{proposition}
Let $\Delta$ be an $n$-dimensional simplicial cone and 
$S=\{\bar\alpha_1, \dots, \bar\alpha_l\}$ be
a set of linear forms on $\Delta$ up to constant multiple such
that $\bar\alpha_i$ is definite on $\Delta$. 
Let $I$ be a definite integral of type $S$.
Then there exists 
\begin{enumerate}
\item
a finite number
of $n$-simplexes $\Delta_1, \dots \Delta_k$,
\item
a flag $F_j$ of $\Delta_j$ for $j=1, \dots, k$,
\item
a derived sequences of linear forms $\Cal D_j =(F_j,\{S^{(i)}_j\})$
($j=1, \dots,k$)
with $S^{(0)}_{j}=S$
for $i=1, \dots, k$ and
\item
and a definite integral $I_j$ of type $\Cal D_j$ for $j=1, \dots,k$
\end{enumerate}
such that
$I$ is $\bold Q^{ab}$-linear combination of $I_1, \dots, I_k$.
\end{proposition}
\begin{proof}
We apply the construction of Proposition \ref{preparation proposition}
for $\Delta$ and 
$S=\{\bar\alpha_1, \dots, \bar\alpha_l\}$
and we get a decomposition 
\begin{equation}
\label{variable change decomposition}
\Delta=\cup_{j \in I}\tilde \Delta_j
\end{equation}
of simplicial cones, a flag $F_j$ on $\tilde\Delta_j$ and a derived
sequence $\Cal D_j=(F_j, \{S^{(i)}_j\})$ with the properties of
Proposition \ref{preparation proposition}. We fix a standard coordinate
$(\eta_1, \dots, \eta_n)$ of $\Delta$. For each $j$, we choose
a coordinate $(\eta_1^{(j)}, \dots, \eta_n^{(j)})$ of $\tilde\Delta_j$
with a relation $\eta_i=\sum_{k=1}^na_{ik}^{(j)}\eta^{(j)}_k$ 
with the property of 
Lemma \ref{changing variable}. We introduce a new variable
$y^{(j)}_1, \dots, y^{(j)}_n$ and define a monomial $y_1, \dots, y_n$
of $y_1^{(j)}, \dots, y_n^{(j)}$ by
\begin{equation}
\label{the relation for exponential}
y_i=y_i(y_1^{(j)},\dots, y_n^{(j)})=
(y_1^{(j)})^{a_{i1}^{(j)}} \cdots (y_n^{(j)})^{a_{in}^{(j)}}
\end{equation}
Since $(0,1)^n=\{\exp(-\eta) \mid \eta \in \Delta\}$,
by the decomposition (\ref{variable change decomposition}),
we have a decomposition of $(0,1)^n=\cup_{j\in I}D_j$, where
\begin{align*}
D_i& =\{\exp(-\eta) \mid \eta \in \tilde\Delta_j \} \\
&=
\{(y_1(y^{(j)}_1,\dots,y^{(j)}_n ), \dots, y_n(y^{(j)}_1,\dots,y^{(j)}_n )
\mid y_1^{(j)} \in (0,1),\dots, y_n^{(j)} \in (0,1)
\}.
\end{align*}
Thus the integral (\ref{definite integral}) is equal to
$$
\int_{(0,1)^n}L(y_1, \dots, y_n)\frac{dy_1}{y_1}\cdots \frac{dy_n}{y_n}
=\sum_{j\in I}
\int_{D_j}L(y_1, \dots, y_n)\frac{dy_1}{y_1}\cdots \frac{dy_n}{y_n}.
$$
For a linear form $\alpha$, a monomial $y^{\alpha}$
can be regarded as a monomial of $y^{(j)}_1, \dots,y^{(j)}_n$ by
the relation (\ref{the relation for exponential}), which is denoted as
$(y^{(j)})^{\alpha^{(j)}}$.
By changing coordinate of integral, we have
\begin{align}
\label{changing variable for integral}
& \int_{D_j}\prod_{i=1}^p
\Big(
\frac{e_iy^{\alpha_i}}{1-e_iy^{\alpha_i}}
\Big)
^{\mu_i}
\frac{dy_1}{y_1}\cdots \frac{dy_n}{y_n} \\
\nonumber
=&
A^{(j)}
\int_{(0,1)^n}\prod_{i=1}^p
\Big(
\frac
{e_i(y^{(j)})^{\alpha_i^{(j)}}}
{1-e_i(y^{(j)})^{\alpha_i^{(j)}}}
\Big)
^{\mu_i}
\frac{dy_1^{(j)}}{y_1}\wedge\cdots \wedge\frac{dy_n^{(j)}}{y_n},
\end{align}
where $A^{(j)}=\det(a^{(j)}_{ij})$.
Then by the property of the coordinate in Lemma \ref{changing variable},
we have
$\alpha^{(j)}_i=e_i{\alpha^{(j)}_i}'$ with
${\alpha^{(j)}_i}'\in S_{j,var}^{(p)}$ for some $p$ and 
$e_i \in \bold N^{\times}$.
Using the equality
$$
\frac{aw^{e}}{1-aw^e}
=\prod_{b^e=a}(1+\frac{bw}{1-bw})-1,
$$
the second line of (\ref{changing variable for integral})
is expressed as a $\bold Q^{ab}$-linear combination
of definite integrals of type $\Cal D_j$.
\end{proof}
\begin{remark}
Shuffle relations for iterated integrals come from these decompositions.
\end{remark}

\section{The third reduction, Definite Integral of type $\Cal D$ and
uni-factor integral of type $\Cal D$}

Let $\Delta$ be an $n$-dimensional simplicial cones, $F$ a flag on
$\Delta$ and $\Cal D=(F, \{S^{(i)}\}_i)$ be a derived sequence of
linear forms.
We consider a function $I(y_n)$ of $y_n$ defined by the following integral.
\begin{equation}
\label{dummy variable}
I(y_n)=
\int_{(0,1)^{n-1}}L(y_1, \dots,y_{n-1},
 y_n) 
\frac{dy_1}{y_1}\cdots \frac{dy_{n-1}}{y_{n-1}}.
\end{equation}
where
$$
L(y_1, \dots, y_n) =\prod_{i=1}^{p}
\frac{e_iy^{\alpha_i}}{(1-e_iy^{\alpha_i})^{\mu_i}} 
\text{with $\alpha_i \in \coprod_{j}S^{(j)}_{var}$, 
$e_i \in \mu_{\infty}$} 
$$
and $\mu_i \geq 1$.
\begin{lemma}
\label{equivalent convergent}
Let $y_n \in (0,1)$. Then the followings are equivalent. 
\begin{enumerate}
\item
The rational function
$$
\frac{L(y_1, \dots ,y_n)}{y_1\cdots y_{n-1}}
$$
is bounded on $(0,1)^{n-1}\times\{y_n\}$.
\item
The numerator of $L(y_1, \dots, y_n)$ is divisible by
$y_i$
for  $1\leq i \leq n-1$.
\item
The integral (\ref{dummy variable}) converges absolutely.

\end{enumerate}

\end{lemma}
The integral (\ref{dummy variable})
is called a uni-factor integral of type $\Cal D$ if
\begin{align*}
L(y_1, \dots, y_n)& =\prod_{i=0}^{n-1}L_i(y_{i+1}, \dots, y_n) \\
L_i(y_{i+1}, \dots, y_n) &=
\begin{cases}
\displaystyle {\frac{e_iy^{\alpha_i}}{(1-e_iy^{\alpha_i})^{\mu_i}} }
\text{with $\alpha_i\in S^{(i)}_{var}$,$\mu_i \geq 1$, 
$e_i \in \mu_{\infty}$, or } \\
1
\end{cases},
\end{align*}
and the equivalent conditions in Lemma \ref{equivalent convergent}
are satisfied.

In the rest of this section, we express the integral
(\ref{dummy variable})
as a linear combination of 
``uni-factor integral''.

\begin{proposition}
The integral $I(y_n)$ in (\ref{dummy variable}) can be written as
a $\bold Q^{ab}$-linear combination of uni-factor integrals of type $\Cal D$.
Moreover if the numerator of $L$ in the expression (\ref{dummy variable})
of the integral $I(y_n)$ is divisible by $y_n$, then
it can be written as a $\bold Q^{ab}$-linear combination of uni-factor
integral such that the numerator of $L$ in (\ref{dummy variable})
is divisible by $y_n$.
\end{proposition}
\begin{proof}
\item
The subset of $\{y_{1}, \dots ,y_{n}\}$ 
consisting of elements dividing the numerator of
$L(y_{1}, \dots , y_{n})$ is called the
zero set of $L(y_{1}, \dots , y_{n})$ and
denoted by $Z(L(y_{1}, \dots, y_{n}))$.
Let $p$ ($1\leq p \leq n$) be an integer.
By the induction of $p$,
we prove that the integral (\ref{dummy variable}) can be expressed as a 
$\bold Q^{ab}$-linear combination of integrals of the form
$$
\int_{(0,1)^{n-p-1}}M^{(p)}(y_{p+1}, \dots, y_n)
\frac{dy_{p+1}}{y_{p+1}}\cdots
\frac{dy_{n-1}}{y_{n-1}}
\int_{(0,1)^p} N^{(p)}(y_1, \dots, y_n)
\frac{dy_{1}}{y_{1}}\cdots
\frac{dy_{p}}{y_{p}},
$$
where
$$
M^{(p)}(y_1, \dots, y_n) =\prod_{i=1}^{l}
\frac{f_iy^{\beta_i}}{(1-f_iy^{\beta_i})^{\nu_i}} 
\text{with $\beta_i \in \coprod_{j\geq p}S^{(j)}_{var}$, 
$\nu_i \geq 1$,$f_i \in \mu_{\infty}$} 
$$ 
\begin{align*}
N^{(p)}(y_1, \dots, y_n)& =\prod_{i=0}^{p-1}N^{(p)}_i(y_{i+1}, \dots, y_n) \\
N^{(p)}_i(y_{i+1}, \dots, y_n) &=
\begin{cases}
\displaystyle {\frac{e_iy^{\alpha_i}}{(1-e_iy^{\alpha_i})^{\mu_i}} }
\text{with $\alpha_i\in S^{(i)}_{var}$, 
$\mu_i \geq 1$,$e_i \in \mu_{\infty}$, or } \\
1
\end{cases},
\end{align*}
with $Z(M^{(p)}N^{(p)}) \supset \{y_1, \dots, y_{n-1}\}$.
To proceed the induction, it is enough to prove the following lemma.
\end{proof}
\begin{lemma}
A rational function
$$
M(y_1, \dots, y_n) =\prod_{i=1}^{l}
\frac{e_iy^{\alpha_i}}{(1-e_iy^{\alpha_i})^{\mu_i}} 
\text{with $\alpha_i \in S^{(p)}_{var}$, 
$e_i \in \mu_{\infty}$} 
$$ 
can be expressed as a $\bold Q^{ab}$-linear combination of 
rational functions of the form
\begin{align*}
N=\frac{fy^{\beta}}{(1-fy^{\beta})^{\nu}} 
\prod_{i=1}^{l}
\frac{e_iy^{\alpha_i}}{(1-e_iy^{\alpha_i})^{\mu_i}} 
\text{with 
$\beta \in S^{(p)}_{var},\alpha_i \in \coprod_{j \geq p+1}S^{(j)}_{var}$, 
$f,e_i \in \mu_{\infty}$} 
\end{align*}
such that $Z(N)\supset Z(M)$.
\end{lemma}
\begin{proof}
We prove the lemma by the induction of $\sum_{i=1}^l \mu_i$.
If $l=1$, then there is nothing to prove. We assume $l\geq 2$ and
$\displaystyle
\frac{e_1y^{\alpha_1}}{(1-e_1y^{\alpha_1})^{\mu_1}} $ and
$\displaystyle
\frac{e_2y^{\alpha_2}}{(1-e_2y^{\alpha_2})^{\mu_2}} $
have distinct factors of denominators.
We put $e_1y^{\alpha_1}=y_{p+1}\gamma_1,e_2y^{\alpha_2}=y_{p+1}\gamma_2$.
By the property of derived sequence, we may assume
$\gamma_2/\gamma_1=(e_2/e_1)\cdot y^{\delta}\in \bold Q[y_{p+2},\dots, y_n]$
and $[\delta]\in S^{(p+1)}$. On the other hand,
we have
\begin{align*}
\frac{\gamma_1 y_{p+1}}{1-\gamma_1 y_{p+1}}  
\frac{\gamma_2 y_{p+1}}{1-\gamma_2 y_{p+1}}=  
\frac{\gamma_2/\gamma_1}{(1-\gamma_2/\gamma_1)}
\frac{\gamma_1 y_{p+1}}{1-\gamma_1 y_{p+1}} 
+
\Big(\frac{-\gamma_2/\gamma_1}{1-\gamma_2/\gamma_1}
-1\Big)
\frac{\gamma_2 y_{p+1}}{1-\gamma_2 y_{p+1}} 
\end{align*}
and
\begin{align*}
Z(\frac{\gamma_1 y_{p+1}}{1-\gamma_1 y_{p+1}}  
\frac{\gamma_2 y_{p+1}}{1-\gamma_2 y_{p+1}}) &=
Z(\frac{\gamma_2/\gamma_1}{(1-\gamma_2/\gamma_1)} 
\frac{\gamma_1 y_{p+1}}{1-\gamma_1 y_{p+1}} ) \\
&=
Z(\frac{\gamma_2 y_{p+1}}{1-\gamma_2 y_{p+1}}).
\end{align*}
By multiplying both sides by
$$
\frac{1}{(1-e_1y^{\alpha_1})^{\mu_1-1}}
\frac{1}{(1-e_2y^{\alpha_2})^{\mu_2-1}}
\prod_{i=1}^{3}
\frac{e_iy^{\alpha_i}}{(1-e_iy^{\alpha_i})^{\mu_i}},
$$
using the relation
$$
\frac{1}{(1-e_iy^{\alpha_i})}=
\frac{e_iy^{\alpha_i}}{(1-e_iy^{\alpha_i})}+1,
$$
we have the statement of the lemma
by the assumption of the induction for $\sum_{i=1}^l \mu_i$.
\end{proof}
\begin{corollary}
For a definite integral $I$ in (\ref{definite integral}),
there exist a $\bold Q^{ab}$-linear combination 
$I(y_n)$ of uni-factor integrals
such that $I(0)=0$ and
$$
I=\lim_{t\to 1}\int_0^t I(y_n)\frac{dy_n}{y_n}.
$$
\end{corollary}

\section{Uni-factor integral of type $\Cal D$ and
Simple uni-factor integral of type $\Cal D$}

Let $\Delta$ be an $n$-dimensional simplex and $F$ be a flag of $\Delta$.
Let $\sigma$ be a regular face of
dimension $n_{\sigma}$ for the flag $F$.
Let $1\leq i_1<\dots <i_{n_\sigma}=n$ be the index set given in
Remark \ref{rem for face and index}
(\ref{equation for sigma}) and 
$(\eta_1^{\sigma},\dots,\eta_1^{\sigma})$ be
the standard coordinate of $\sigma$.
We consider an embedding
$$
(\bold C^{\times})^{n_{\sigma}} \subset (\bold C^{\times})^n
$$
defined by
$$
(y_1^{\sigma}, \dots, y_{n_{\sigma}}^{\sigma})
\mapsto (1, \dots, 1,\overset{\overset{i_1}\smile}
{ y_{1}^{\sigma}}, 1,\dots ,1,
\overset{\overset{i_2}\smile}
{ y_{2}^{\sigma}}, 1,\dots,1, 
\overset{\overset{i_{n_\sigma}}\smile}{ y_{n_\sigma}^{\sigma}})
$$ 
We generalize the notion of uni-factor integral by introducing
the notion of weight.
Let $\Cal D$ be a derived sequence of linear forms and the restriction
of $\Cal D$ to $\sigma$ is denoted by 
$\Cal D_{\sigma}=(F,\{S^{(i)}_{\sigma}\})$.
Let $k \leq n-1$.
We consider a rational function
\begin{align}
L=L(y_1^{\sigma}, \dots, y_{n_{\sigma}}^{\sigma})& 
=\prod_{i=0}^{n_\sigma-1}L_i(y_{i+1}^{\sigma}, \dots, y_{n_\sigma}^{\sigma}) \\
\label{def of unifactor 2}
L_i(y_{i+1}^{\sigma}, \dots, y_{n_{\sigma}}^{\sigma}) &=
\begin{cases}
\displaystyle {\frac{e_i{y^{\sigma}}^{\alpha_i}}
{(1-e_i{y^{\sigma}}^{\alpha_i})^{\mu_i}} }
\text{with $\alpha_i\in S^{(i)}_{\sigma,var}$, 
$e_i \in \mu_{\infty}$, or } \\
1
\end{cases},
\end{align}
satisfying the condition 
\begin{equation}
\label{sigma bounded condition}
\text{The numerator of $L(y_1^{\sigma}, \dots, y_{n_\sigma}^{\sigma})$ 
is divisible by
$y_i^{\sigma}$}
\text{ for } 1\leq i \leq k.
\end{equation}
Since $\prod_{i=k}^{n_\sigma-1}L_i$ is independent of 
$y_1^{\sigma},\dots,y_k^{\sigma}$, the condition 
(\ref{sigma bounded condition}) is equivalent to the condition
\begin{equation*}
\text{The numerator of 
$\prod_{i=1}^{k-1}L_i(y_{i+1}^{\sigma}, \dots, y_{n_{\sigma}}^{\sigma})$ 
is divisible by
$y_i^{\sigma}$}
\text{ for } 1\leq i \leq k.
\end{equation*}
Then the integral 
\begin{equation}
\label{sigma dummy variable}
I(y_{k+1}^{\sigma},y_{k+2}^{\sigma},\dots, y_{n_{\sigma}}^{\sigma})
=
\int_{(0,1)^{k}}L(y_1^{\sigma}, \dots,
 y_{n_{\sigma}}^{\sigma}) 
\frac{dy_1^{\sigma}}{y_1^{\sigma}}\cdots 
\frac{dy_{k}^{\sigma}}{y_{k}^{\sigma}}.
\end{equation}
is a function of 
$
(y_{k+1}^{\sigma},y_{k+2}^{\sigma},\dots, y_{n_{\sigma}}^{\sigma}).
$
\begin{definition}
\begin{enumerate}
\item
The integral (\ref{sigma dummy variable}) 
is called a uni-factor integral of type $\Cal D_{\sigma}$ of weight $k$.
\item
A uni-factor integral is said to be simple,
if $L_i=1$ or $\mu_i=1$ for $i\leq k-1$
and $L_k=L_{k+1}=\cdots = L_{n_\sigma-1}=1$
in the expression (\ref{def of unifactor 2}). 
\item
The subset of $\{y_{1}^{\sigma}, \dots ,y_{n_\sigma}^{\sigma}\}$ 
consisting of elements dividing the numerator of
$L(y_{1}^{\sigma}, \dots , y_{n_\sigma}^{\sigma})$ is called the
zero set of $L(y_{1}^{\sigma}, \dots , y_{n_\sigma}^{\sigma})$ and
denoted by $Z(L(y_{1}^{\sigma}, \dots, y_{n_\sigma}^{\sigma}))$.
For a uni-factor integral (\ref{sigma dummy variable}),
$Z(L(y_{1}^{\sigma}, \dots, y_{n_\sigma}^{\sigma}))
-\{y_{1}^{\sigma}, \dots, y_{k_\sigma}^{\sigma}\}$ is 
called the zero set of the uni-factor integral and denoted 
by
$Z(I(y_{k+1}^{\sigma}, \dots, k_{n_\sigma}^{\sigma}))$.
\end{enumerate}
\end{definition}
If $y_i\in Z(I(y_{k+1}^{\sigma}, \dots, k_{n_\sigma}^{\sigma}))$,
then
$
I(y_{k+1}^{\sigma}, \dots, k_{n_\sigma}^{\sigma})\mid_{y_i=0}=0.
$
The set of functions on $(y_{k+1}^{\sigma}, \dots, y_{n_\sigma}^{\sigma})$ 
generated by
uni-factor (resp. simple uni-factor ) integrals of
weight $k$ is
denoted by $\Cal U_{\sigma,k}$
(resp. $\Cal S_{\sigma,k}$).

We consider three sequence of statements.
We use $n=n^{\sigma}$ and $(y_1, \dots, y_n)$
for $(y_1^{\sigma}, \dots ,y_n^{\sigma})$ for simplicity.
\begin{enumerate}
\item [($A_{\sigma,k}$)]
Let $h(y_{k+1}, \dots ,y_n)$ be an element of $\Cal U_{\sigma,k}$
and set $Z=Z(h(y_{k+1}, \dots ,y_n))$.
Then $h(y_{k+1}, \dots ,y_n)$ can be written as a 
$\bold Q^{ab}$-linear combination of
\linebreak
$M(y_{k+1}, \dots,y_n)f(y_{k+1}, \dots, y_n)$, where
\begin{enumerate}
\item
$f(y_{k+1}, \dots, y_n) \in \Cal S_{\sigma',k'}$ with 
$k' \leq k, \sigma' \leq \sigma$,
\item
\begin{align}
M(y_{k+1}, \dots,y_n)& =\prod_{i=k}^{n-1}
M_i(y_{i+1}, \dots,y_n), \\
\label{each factor}
M_i(y_{i+1}, \dots,y_n)& 
=
\begin{cases}
\displaystyle
\frac{ey^\alpha_i}{(1-ey^\alpha_i)^{\nu_{\alpha_i}}} 
\text{ with $\alpha_i \in S^{(i)}_{\sigma,var}$ }
\\
1,
\end{cases}
\end{align}
and
\item
$Z \subset Z(M(y_{k+1}, \dots,y_n)f(y_{k+1}, \dots, y_n))$
\end{enumerate}

\item[($B_{\sigma,k}$)]
Let $h(y_{k+1}, \dots, y_n)$ be an element of $\Cal S_{\sigma,k}$
and set $Z=Z(h(y_{k+1}, \dots, y_n))$.
Then for any $i\in [k+1,n]$,
$\displaystyle y_i\frac{\partial}{\partial y_i}h(y_{k+1}, \dots, y_n)$
can be written as a $\bold Q^{ab}$-linear combination of 
$f(y_{k+1}, \dots, y_n) \in \Cal U_{\sigma',k'}$
with $k'<k$
and
$Z \subset Z(f(y_{k+1}, \dots, y_n))$

\item[($C_{\sigma}$)]
The statement ($A_{\sigma,k}$) and ($B_{\sigma,k}$) holds for all $k$.
\end{enumerate}
\begin{proposition}
\label{induction A}
$$ 
\begin{matrix}
(C_{\sigma'} \text{ for } \sigma' < \sigma) 
\text{ and } \\
(A_{\sigma,k'} \text{ for } k' < k) 
\text{ and } \\
(B_{\sigma,k'} \text{ for } k' < k) 
\end{matrix}
\Rightarrow
(A_{\sigma,k}) 
$$
\end{proposition}
\begin{proof}
Let $f=f(y_{k+1},\dots, y_n)$ be an element of $\Cal U_{\sigma,k}$.
Then $f$ can be written as
$$
f(y_{k+1},\dots, y_n)=\int_0^1g(y_k, \dots, g_n)
\frac{dy_k}{y_k},
$$
where $g \in \Cal U_{\sigma, k-1}$.
We put
$Z(g(y_k, \dots, y_n))=Z\cup \{y_k\}$.
By the inductive hypothesis, the function $g(y_k, \dots,y_1)$
can be written as a $\bold Q^{ab}$-linear combination of
$
M(y_k,\dots, y_n)h(y_k, \dots, y_n)
$
where (1) $h(y_k, \dots, y_n) \in \Cal S_{\sigma, var}^{(k-1)}$,
(2)$M=\prod_{i=k-1}^{n-1}M_i$
where $M_i$ is given in the form (\ref{each factor}), and (3)
$
Z(g(y_k, \dots, y_n)) \subset
Z(M(y_k,\dots, y_n)h(y_k, \dots, y_n))$.

\noindent
(I)
The case $M_{k-1}=1$. In this case $y_k \in Z(h)$ and
the integral
$$
\int_{0}^1h(y_k, \dots, y_n)\frac{dy_k}{y_k}
$$
is simple uni-factor integral by definition.

\noindent
(II)
The case where
the $(k-1)$-th factor of $M$ is given as
$M_{k-1}=\displaystyle \frac{ey^{\alpha_{k-1}}}
{(1-ey^{\alpha_{k-1}})^{e_{k-1}}}$.
We put $y^{\alpha_{k-1}}=py_k$ with a monomial
$p=p(y_{k+1},\dots, y_n)$.
If $e_{k-1}=1$,
the integral
\begin{equation}
\label{reduction A e=1}
\int_0^1\frac{epy_k}
 {(1-epy_k)^{e_{k-1}}}
h(y_k, \dots, y_n)\frac{dy_k}{y_k} 
\end{equation}
is an element in $\Cal S_{\sigma,k}$. If $e_{k-1} >1$, then the integral
(\ref{reduction A e=1})
is equal to
\begin{align}
\label{partial int1}
&
\frac{1}{e_{k-1}-1}\Big[ \big(
\frac{1}{(1-epy_k)^{e_{k-1}-1}}-1 \big) h(y_k, \dots, y_n)\Big]_{y_k=0}^1 \\
\label{partial int2}
-&
\frac{1}{e_{k-1}-1}\int_0^1 \big(
\frac{1}{(1-epy_k)^{e_{k-1}-1}}-1 \big) 
\big( y_k
\frac{\partial}{\partial y_k} h(y_k, \dots, y_n)\big)
\frac{dy_k}{y_k}
\end{align}
The first term (\ref{partial int1}) can be written
as a uni-factor integral of weight less than $k$ by 
the equality
$$
\frac{1}{(1-x)^m}=
1+\frac{x}{1-x}+\frac{x}{(1-x)^2}+\cdots +\frac{x}{(1-x)^m}.
$$
By the assumption ($B_{\sigma, k-1}$),
the second term (\ref{partial int2})
can be written as a $\bold Q^{ab}$-linear combination of
uni-factor integral of weight less than $k$.
We can see that the zero set of each term in
(\ref{partial int1}) and (\ref{partial int2})
contains $Z(M_{k-1}h)-\{y_k\}$.
\end{proof}

\begin{proposition}
\label{induction B}
$$
(B_{\sigma,k'} \text{ for } k' < k) \Rightarrow
(B_{\sigma,k}) 
$$
\end{proposition}
\begin{proof}
Let $f=f(y_{k+1},\dots, y_n)$ be an element of $\Cal S_{\sigma,k}$.
Then $f$ can be written as
$$
f(y_{k+1},\dots, y_n)=\int_0^1L_{k-1}(y_k, \dots,y_n)g(y_k, \dots, g_n)
\frac{dy_k}{y_k},
$$
where $g \in \Cal S_{\sigma, k-1}$,
$\displaystyle L_{k-1}=\frac{ey^{\alpha_{k-1}}}
{1-ey^{\alpha_{k-1}}}$
with $\alpha_{k-1} \in S_{\sigma}^{(k-1)}$ or $1$.

\noindent
(I) If $L_{k-1}=1$, then
$Z(f)=Z(g)-\{y_k\}$ and 
$$
y_i\frac{\partial}{\partial y_i}
f(y_{k+1},\dots, y_n)=
\int_0^1y_i\frac{\partial}{\partial y_i}g(y_k, \dots, g_n)
\frac{dy_k}{y_k},
$$
and by inductive hypothesis,
it is a $\bold Q^{ab}$-linear combination of uni-factor integral 
of weight less than $k$.

\noindent
(II) If $\displaystyle L_{k-1}=\frac{ey^{\alpha_{k-1}}}
{1-ey^{\alpha_{k-1}}}$, then
$Z(f)=Z(L_{k-1}g(y_k, \dots, y_1))$. We have
\begin{align}
\nonumber
 y_i\frac{\partial}{\partial y_i}
f(y_{k+1},\dots, y_n) 
& =
\int_0^1
\frac{ey^{\alpha_{k-1}}}
{1-ey^{\alpha_{k-1}}}
y_i\frac{\partial}{\partial y_i}g(y_k, \dots, g_n)
\frac{dy_k}{y_k} \\
\label{simple term 1}
& +\int_0^1
\frac{cey^{\alpha_{k-1}}}
{(1-ey^{\alpha_{k-1}})^2}
g(y_k, \dots, g_n)
\frac{dy_k}{y_k} 
\end{align}
Here $c$ is the coefficient of $\alpha_{k-1}$ on $\eta_i$.
The second term (\ref{simple term 1}) is equal to 
$$
\Big[
\frac{cey^{\alpha_{k-1}}}
{1-ey^{\alpha_{k-1}}}
g(y_k, \dots, g_n)
\Big]_{y_k=0}^1
-\int_0^1
\frac{cey^{\alpha_{k-1}}}
{1-ey^{\alpha_{k-1}}}
y_k\frac{\partial}{\partial y_k}g(y_k, \dots, g_n)
\frac{dy_k}{y_k}. 
$$
By inductive hypothesis, we have the proposition.
\end{proof}

\section{Proof of the main theorem}

By Proposition \ref{induction A}
and Proposition \ref{induction B}, we have the following proposition
\begin{proposition}
\label{functional recursion cone}
\begin{enumerate}
\item
An element $f(y_n)$ in $\Cal U_{\Delta,n-1}$ can be written
as a $\bold Q^{ab}$-linear combination of
$$\frac{1}{(1-ey_n)^\mu}g(y_n)$$
with $\mu \geq 0$, $g(y_n)\in \Cal S_{\Delta,n-1}$.
\item
For an element $f(y_n) \in \Cal S_{\Delta,n-1}$,
$$
y_n\frac{\partial}{\partial y_n}f(y_n)
$$
is an element of $\Cal U_{\Delta,n-2}$.
\end{enumerate}
\end{proposition}
\begin{definition}
\begin{enumerate}
\item
For a differential form $\omega_1, \dots, \omega_k$ on $x$,
an iterated integral is defined inductively by
$$
\int_0^y\omega_1\omega_2\cdots\omega_k
=\int_0^y(\omega_1(z)\int_0^z\omega_2\cdots\omega_k)
$$
if $k\geq 2$ and usual one if $k=1$. 
\item
A function on $y$
defined by
$$
\int_0^y
(\frac{dx}{x})^{k_1-1}\frac{dx}{1-e_1x}
(\frac{dx}{x})^{k_2-1}\frac{dx}{1-e_2x}\cdots
(\frac{dx}{x})^{k_m-1}\frac{dx}{1-e_mx},
$$
where $k_i \geq 1$, $e\in \mu_{\infty}$ is called a multiple polylogarithm
of weight $k=k_1+\cdots + k_m$.
\item
Let $\Cal P_k$ be a $\bold Q^{ab}$-linear combination of
$\displaystyle\frac{1}{(1-ey)^{\mu}}h(y)$, 
where $\mu \geq 0$ and $h(y)$ is a multiple 
polylogarithm of weight less than or equal
to $k$.
\end{enumerate}
\end{definition}
\begin{proposition}
\label{functional recursion multiple polylog}
Let $f(y)$ be an element of $\Cal P_k$. Then
$$
\int_{0}^y
\frac{1}{(1-et)^{\nu}}f(t)dt,\quad (\nu \geq 1)\qquad
\int_{0}^y
\frac{1}{t}f(t)dt
$$
are elements of $\Cal P_{k+1}$.
\end{proposition}
\begin{proof}
We prove by the induction of $k$.
Let $f(y)$ be a multiple polylogarithm of weight $k$
and we show that the integral
$
\displaystyle
\int_{0}^y\frac{1}{(1-et)^{\nu}}f(t)dt
$
is an element of $\Cal P_{k+1}$ by the induction on $\nu$.
If $\nu=1$, the statement is true by the definition of
multiple polylogarithm.
If $\nu \geq 2$, we have
\begin{align}
\label{zero multiple polylog}
&\int_{0}^y\frac{1}{(1-et)^{\nu}}f(t)dt \\
\label{first multiple polylog}
& =
\big[
\frac{1}{e(\nu-1)(1-et)^{\nu-1}}f(t)
\big]_0^y
\\
\label{second multiple polylog}
& -
\frac{1}{e(\nu-1)}
\int_{0}^y\frac{1}{(1-et)^{\nu-1}}\frac{d}{dt}f(t)dt
\end{align}
The first term (\ref{first multiple polylog}) 
is an element in $\Cal P_{k+1}$ by the hypothesis of induction
on $\nu$.
The second term
(\ref{second multiple polylog}) is an element in $\Cal P_k$ by 
the hypothesis of induction on $k$. As a consequence, the integral
(\ref{zero multiple polylog})
is an element of $\Cal P_{k+1}$.
\end{proof}
By Proposition \ref{functional recursion cone}, 
Proposition \ref{functional recursion multiple polylog}, 
we have the following theorem.
\begin{theorem}
\label{definite and polylog}
$\Cal U_{\Delta,k} \subset\Cal P_k$.
\end{theorem}
\begin{theorem}
If the definite integral (\ref{definite integral}) 
is absolutely convergent,
it is an element in $\Cal Z_{\infty}$. 
\end{theorem}
\begin{proof}
We apply Theorem \ref{definite and polylog} to the function
$I(y_n)$ defined in (\ref{dummy variable}).
We put 
$$
f(y)=\int_{0}^yI(y_n)\frac{dy_n}{y_n}
$$
Then by Proposition \ref{functional recursion multiple polylog},
$f(y)\in \Cal P_{k}$ for some $k$ and 
$\displaystyle \lim_{y\to 1}f(y)=I$
exists.
Since $f(y)$ can be written as
$$
\sum_{m<0}\frac{1}{(1-y)^m}h_m(y) +
\sum_{\{(e,m)\mid e \neq 1\text{ or }m=0\}}\frac{1}{(1-ey)^m}H_{e,m}(y)
$$
where $h_m(y)$ and $H_{e,m}(y)$ is a $\bold Q^{ab}$-linear
combination of multiple polylogarithm.
By the theory of regularization in \cite{IKZ}, \cite{R},
$h_m(y)$ and $H_{e,m}$ can be written as
\begin{align*}
h_m(y) &=\sum_{i=0}^{N}a_{m,i}(-\log (1-y))^i+0((1-x)^{\epsilon}) \\
H_{e,m}(y) &=\sum_{i=0}^{N}a_{e,m,i}(-\log (1-y))^i+0((1-x)^{\epsilon}),
\end{align*}
where $a_{m,i}$ and $a_{e,m,i}$ is a $\bold Q^{ab}$-linear
combination of cyclotomic multiple zeta values.
(\cite{IKZ},\cite{R}).
Since $\displaystyle\lim_{y\to 1}f(y)$ exists,
we have $a_{m,i}=0$ for $m >0, i\geq 0$ and $a_{e,m,i}=0$ for $i>0$.
Therefore we have 
$$
I=\sum_{\{(e,m)\mid e \neq 1\text{ or }m=0\}}
\frac{a_{e,m,0}}{(1-e)^m}.
$$
\end{proof}
 
\end{document}